\documentclass[a4paper,11pt,twoside,reqno]{amsart}

\usepackage[utf8]{inputenc}
\usepackage[plainpages=false,pdfpagelabels=true]{hyperref}
\usepackage{amssymb,amsthm}
\usepackage[margin=1in]{geometry}

\newtheorem{Satz}{Theorem}[section]
\newtheorem{Prop}[Satz]{Proposition}
\newtheorem{Lem}[Satz]{Lemma}

\newtheorem{Cor}[Satz]{Corollary}
\theoremstyle{definition}

\newtheorem{Bem}[Satz]{Remark}
\newcommand{\zarg}{d\phi(\partial_t)\wedge d\phi(\partial_x)}
\newcommand{\Hom}{\operatorname{Hom}}
\newcommand{\hess}{\operatorname{Hess}}
\newcommand{\tr}{\operatorname{Tr}}
\newcommand{\sff}{\mathrm{I\!I}}

\allowdisplaybreaks[1]

\renewcommand{\epsilon}{\varepsilon}

\newcommand{\R}{\ensuremath{\mathbb{R}}}

\numberwithin{equation}{section}

\title{On the full bosonic string from Minkowski space to Riemannian manifolds}
\author{Volker Branding}
\date{\today}
\address{University of Vienna, Faculty of Mathematics\\
Oskar-Morgenstern-Platz 1, 1090 Vienna, Austria}
\email[]{volker.branding@univie.ac.at}
\subjclass[2010]{35L71; 58J45; 35L05}
\keywords{wave maps with scalar and two-form potential, full bosonic string, Cauchy problem}
\begin{document}

\begin{abstract}
We study the action of the full bosonic string for the domain being two-dimensional Minkowski space
and the target a Riemannian manifold.
Its critical points couple the wave map equation to a scalar and a two-form potential.
Besides investigating their basic features we establish existence results for the latter.
\end{abstract} 

\maketitle

\section{Introduction and Results}
The nonlinear \(\sigma\)-model is one of the most prominent models studied in quantum field theory.
In mathematical terms it is described by the Dirichlet energy of a map
between two manifolds. In the case that the domain has a Riemannian metric its critical points
lead to the \emph{harmonic map equation}, which is a semilinear elliptic equation,
see \cite{MR2389639} for a recent survey.
However, if the domain has a Lorentzian metric its critical points are given by the \emph{wave map equation},
which is a semilinear hyperbolic equation, see \cite{MR1674843,MR2043751,MR3098646} and references therein
for an introduction to wave maps.

In this article we will focus on the nonlinear \(\sigma\)-model coupled to two external potentials,
which arises as the full action of the bosonic string, see \cite{MR2151029}, p. 108 for the physical background.
More precisely, we will study the following action
\begin{equation}
\label{energy-functional}
S(\phi)=\int_M\big(\frac{1}{2}|d\phi|^2+\phi^\ast B+V(\phi)\big)d\mu.
\end{equation}
Here, \(\phi\colon M\to N\) is a map between two manifolds \(M\) and \(N\) with \(\dim M=2\) and \(\dim N\geq 3\).
In addition, \(B\) is a two-form on \(N\), which we pullback via the map \(\phi\).
Finally, \(V\colon N\to\R\) is a scalar function. For the sake of completeness, we want to mention that 
in the physics literature one often finds that the scalar potential \(V(\phi)\) in \eqref{energy-functional} gets multiplied
by the scalar curvature of the domain. The action functional \eqref{energy-functional} models a bosonic string under the influence
of an external magnetic and scalar field.

This article is a sequel to \cite{MR3573990}, where the action functional \eqref{energy-functional} has been investigated
in the case that the domain is a closed Riemannian surface and the target a Riemannian manifold.
In this case an existence result could be obtained by the heat-flow method, which requires the target manifold to have
negative curvature.

Here, we consider the situation that the domain \(M\) is two-dimensional Minkowski space and the target a  
Riemannian manifold. It turns out that contrary to the Riemannian case
we do not have to impose a curvature condition on the target to obtain an existence result.

For wave maps from two-dimensional Minkowski space to Riemannian manifolds there exist two classical
existence results. The first one is due to Gu \cite{MR596432}, who showed
that the Cauchy problem for wave maps from two-dimensional Minkowski space has a global smooth solution
for smooth initial data. Using the method of a priori estimates the existence of a global weak solution
for special target manifolds was shown in \cite{MR678488}, which was later extended to arbitrary targets 
in \cite{MR853598}.

A problem similar to the one studied in this article has been investigated in \cite{MR2334944,MR1303178}:
The authors study timelike minimal surfaces in Lorentzian manifolds, which corresponds to
critical points of \eqref{energy-functional} without potentials and the target being a Lorentzian manifold.

This article is organized as follows: In Section 2 we study the basics of the action functional \eqref{energy-functional}.
In Section 3 we derive an existence result for critical points of \eqref{energy-functional} in the special case of a vanishing 
scalar potential \(V(\phi)\) making use of the method of characteristics. In addition, we point out how this result can be obtained
by purely analytical methods.
In the last section we prove the existence of a global weak solution to \eqref{energy-functional} for a non-vanishing scalar potential \(V(\phi)\).

\section{Wave maps with scalar and two-form potential}
In this section we investigate the basic features of \eqref{energy-functional}.

Let \((M,h)\) be two-dimensional Minkowski space with global coordinates \((t,x)\) and the metric of signature \((+,-)\). 
Moreover, let \((N,g)\) be a Riemannian manifold of dimension \(\dim N=n\geq 3\).
We will use Latin letters for indices on the target \(N\) and Greek letters for indices on the domain \(M\).
In addition, we make use of the summation convention, that is we sum over repeated indices.

\begin{Prop}
The Euler-Lagrange equation of the functional \eqref{energy-functional} is given by
\begin{equation}
\label{euler-lagrange}
\tau(\phi)=Z(\zarg)+\nabla V(\phi),
\end{equation}
where \(\tau(\phi):=\tr(\nabla^{\phi^\ast TN}d\phi)\) denotes the wave map operator of the map \(\phi\). The vector-bundle homomorphism \(Z \in \Gamma (\Hom(\Lambda^2T^\ast N,TN))\) is defined by the equation
\begin{equation}
\label{def-Z}
\Omega(\eta,\xi_1,\xi_2)=\langle Z(\xi_1\wedge\xi_2),\eta\rangle,
\end{equation}
where \(\Omega=dB\) is a three-form on \(N\) and \(\eta,\xi_1,\xi_2\in TN\).
Moreover \(\partial_t,\partial_x\) is the global orthonormal basis of \(T\R^{1,1}\simeq\R^{1,1}\).
\end{Prop}
\begin{proof}
This follows as in the Riemannian case, see \cite[Proposition 2.1]{MR3573990}.
\end{proof}

We call solutions of \eqref{euler-lagrange} \emph{wave maps with scalar and two-form potential}.

In terms of global coordinates \((t,x)\) on two-dimensional Minkowski space and local coordinates \(y^i,i=1\ldots,n\) on the target \(N\)
the Euler-Lagrange equation \eqref{euler-lagrange} acquires the form

\begin{equation}
\label{wavemap-b-v-local}
\frac{\partial^2\phi^i}{\partial t^2}-\frac{\partial^2\phi^i}{\partial x^2}
+\Gamma^i_{jk}(\frac{\partial\phi^j}{\partial t}\frac{\partial\phi^k}{\partial t}-\frac{\partial\phi^j}{\partial x}\frac{\partial\phi^k}{\partial x})
-Z^i(\partial_{y^k}\wedge\partial_{y^j})\frac{\partial\phi^k}{\partial t}\frac{\partial\phi^j}{\partial x}-g^{ij}\frac{\partial V(\phi)}{\partial y^j}=0.
\end{equation}
Here, \(\Gamma^i_{jk}\) are the Christoffel symbols of the Levi-Cevita connection on \(N\).

\begin{Bem}
The standard wave map equation can also be studied from higher-dimensional Minkowski space.
However, the two-form potential in \eqref{energy-functional} only gives a well-defined action functional
in two dimensions such that the study of \eqref{energy-functional} is restricted to a two-dimensional domain.
\end{Bem}

\begin{Bem}
The vector-bundle homomorphism \(Z\) defined in \eqref{def-Z} can be interpreted as arising from a metric connection with totally antisymmetric torsion.
In this case one has
\[
\nabla^{Tor}_XY=\nabla^{LC}_XY+A(X,Y),
\]
where \(\nabla^{LC}\) denotes the Levi-Cevita connection, \(X,Y\) are vector fields and \(A(X,Y)\) is a skew-adjoint endomorphism.
The endomorphism \(A(X,Y)\) satisfies
\[
\langle A(X,Y),Z\rangle=\Omega(X,Y,Z)
\]
with \(\Omega\in\Gamma(\Lambda^3T^\ast N)\) similar to \eqref{def-Z}. For more details on metric connections with torsion
see \cite{MR3493217} and references therein.
\end{Bem}

By varying the action functional with respect to the domain metric we obtain the energy-momentum tensor:
\begin{equation*}
T_{\alpha\beta}:=\frac{1}{2}h_{\alpha\beta}|d\phi|^2-\langle d\phi(e_\alpha),d\phi(e_\beta)\rangle+V(\phi)h_{\alpha\beta}.
\end{equation*}

By a direct calculation it follows that the energy-momentum tensor \(T_{\alpha\beta}\) is divergence free whenever \(\phi\)
is a wave map with scalar and two-form potential, see \cite[Lemma 3.3]{MR3573990} for more details.

For many of the calculations in this article it turns out to be useful to exploit the geometry of two-dimensional Minkowski space 
by introducing lightcone coordinates:
\[
\xi=\frac{x+t}{2},\qquad \eta=\frac{x-t}{2},\qquad \frac{\partial}{\partial x}=\frac{1}{2}\big(\frac{\partial}{\partial\xi}+\frac{\partial}{\partial\eta}\big)
,\qquad \frac{\partial}{\partial t}=\frac{1}{2}\big(\frac{\partial}{\partial\xi}-\frac{\partial}{\partial\eta}\big).
\]

In these coordinates \eqref{wavemap-b-v-local} acquires the form
\begin{align}
\frac{\partial^2\phi^i}{\partial\xi\partial\eta}+(\Gamma^i_{jk}
+\frac{1}{2}Z^i(\partial_{y^k}\wedge\partial_{y^j}))\frac{\partial\phi^j}{\partial\xi}\frac{\partial\phi^k}{\partial\eta}+g^{ij}\frac{\partial V(\phi)}{\partial y^j}=0,\qquad i=1,\ldots,n.
\end{align}

At several places it becomes necessary to make use of the Nash embedding theorem to isometrically embed \(N\) into some \(\R^q\).

\begin{Lem}
Assume that \(N\subset\R^q\). Then the Euler-Lagrange equation \eqref{euler-lagrange} acquires the form
\begin{equation}
\label{euler-lagrange-extrinsic}
\frac{\partial^2u}{\partial t^2}-\frac{\partial^2u}{\partial x^2}=\sff(du,du)+\tilde{Z}(du(\partial_t)\wedge du(\partial_x))+\widetilde{\nabla V}(u),
\end{equation}
where \(u\colon\R^{1,1}\to\R^q\) and \(\sff\) denotes the second fundamental form of \(N\) in \(\R^q\). Moreover, \(\tilde{Z}\) and \(\widetilde{\nabla V}\) 
denote the extensions of \(Z\) and \(\nabla V\) to the ambient space \(\R^q\).
In terms of lightcone coordinates on \(M\) the Euler-Lagrange equation \eqref{euler-lagrange} is given by
\begin{align}
\label{euler-lagrange-extrinsic-lightcone}
\frac{\partial^2u}{\partial\xi\partial\eta}=\sff(\frac{\partial u}{\partial\eta},\frac{\partial u}{\partial\xi})-\frac{1}{2}\tilde{Z}(du(\partial_\xi)\wedge du(\partial_\eta))
-\widetilde{\nabla V}(u).
\end{align}
\end{Lem}

\begin{proof}
For a proof see \cite[Lemma 3.8]{MR3573990} and references therein.
\end{proof}

In the following we will omit the tildes at \(Z\) and \(\nabla V\) when extending them to the ambient space.

Besides the energy-momentum tensor there is another conserved quantity:
\begin{Prop}
Let \(u\colon\R^{1,1}\to\R^q\) be a solution of \eqref{euler-lagrange-extrinsic}.
Then the energy
\begin{align}
\label{conserved-energy}
E(t)=\int_\R(\big|\frac{\partial u}{\partial t}\big|^2+\big|\frac{\partial u}{\partial x}\big|^2-2V(u))dx
\end{align}
is constant with respect to \(t\).
\end{Prop}
\begin{proof}
Using the extrinsic version of the equation for wave maps with scalar and two-form potential \eqref{euler-lagrange-extrinsic}
we calculate
\begin{align*}
\frac{d}{dt}\int_\R(\big|\frac{\partial u}{\partial t}\big|^2+\big|\frac{\partial u}{\partial x}\big|^2)dx
=&2\int_\R\langle\frac{\partial u}{\partial t},\frac{\partial^2 u}{\partial t^2}-\frac{\partial^2 u}{\partial x^2}\rangle dx \\
=&2\int_\R\big(\underbrace{\langle\frac{\partial u}{\partial t},\sff(du,du)\rangle}_{=0\text{ since } \sff\perp\partial_u}
+\langle\frac{\partial u}{\partial t},Z(\frac{\partial u}{\partial t}\wedge\frac{\partial u}{\partial x})\rangle 
+\langle\frac{\partial u}{\partial t},\nabla V(u)\rangle\big)dx \\
=&2\int_\R\big(\underbrace{\Omega(\frac{\partial u}{\partial t},\frac{\partial u}{\partial t},\frac{\partial u}{\partial x})}_{=0}+\frac{\partial}{\partial t}V(u)\big)dx,
\end{align*}
where we interchanged the derivatives with respect to \(t\) and \(x\) and used integration by parts in the first step.
\end{proof}

\section{Wave maps with two-form potential}
In this section we consider the simpler case of a vanishing scalar potential, that is \(V=0\).
The action functional under consideration then reduces to
\[
S(\phi)=\int_{\R^{1,1}}(\frac{1}{2}|d\phi|^2+\phi^\ast B)d\mu
\]
and its critical points are given by
\begin{align}
\label{euler-lagrange-magnetic}
\tau(\phi)=Z(\zarg). 
\end{align}
We call solutions of \eqref{euler-lagrange-magnetic} \emph{wave maps with two-form potential}.

Note that solutions of \eqref{euler-lagrange-magnetic} are invariant under scaling in the sense that
if \(\phi(t,x)\) is a solution of \eqref{euler-lagrange-magnetic}, then \(\phi_\lambda(t,x)\) defined by
\[
\phi_\lambda(t,x):=\phi(\lambda t,\lambda x)
\]
is also a solution of \eqref{euler-lagrange-magnetic}, where \(\lambda\) is a positive real-valued parameter.

\begin{Bem}
In the case that the domain is a Riemannian surface, the target a three-dimensional
Riemannian manifold and \(\phi\) an isometric embedding, 
the equation \eqref{euler-lagrange-magnetic} is known as \emph{prescribed curvature equation}.
Since we are considering a domain with an indefinite metric here, it does not seem to be possible
to give an analogue geometric interpretation of solutions of \eqref{euler-lagrange-magnetic}.
\end{Bem}

By Proposition \ref{conserved-energy}, using the extrinsic version of \eqref{euler-lagrange-magnetic}, the energy 
\[
E(t):=\frac{1}{2}\int_\R(\big|\frac{\partial\phi}{\partial t}\big|^2+\big|\frac{\partial\phi}{\partial x}\big|^2)dx
\]
is constant with respect to \(t\) whenever \(\phi\) is a solution of \eqref{euler-lagrange-magnetic}.

Locally, \eqref{euler-lagrange-magnetic} looks like
\begin{equation}
\label{wavemap-b-local}
\frac{\partial^2\phi^i}{\partial t^2}-\frac{\partial^2\phi^i}{\partial x^2}
+\Gamma^i_{jk}(\frac{\partial\phi^j}{\partial t}\frac{\partial\phi^k}{\partial t}-\frac{\partial\phi^j}{\partial x}\frac{\partial\phi^k}{\partial x})
-Z^i(\partial_{y^k}\wedge\partial_{y^j})\frac{\partial\phi^k}{\partial t}\frac{\partial\phi^j}{\partial x}=0,\qquad i=1,\ldots,n
\end{equation}
and in terms of lightcone coordinates \((\xi,\eta)\) equation \eqref{wavemap-b-local} acquires the form
\begin{align}
\frac{\partial^2\phi^i}{\partial\xi\partial\eta}+(\Gamma^i_{jk}+\frac{1}{2}Z^i(\partial_{y^k}\wedge\partial_{y^j}))\frac{\partial\phi^j}{\partial\xi}\frac{\partial\phi^k}{\partial\eta}=0,\qquad i=1,\ldots,n.
\end{align}
We will establish the following existence result:

\begin{Satz}
\label{satz-cauchy-magnetic}
Let \((N,g)\) be a complete Riemannian manifold. Then the Cauchy problem 
for wave maps with two-form potential, that is
\begin{align}
\label{cauchy-magnetic}
\tau(\phi)=Z(\zarg)
\end{align}
with the smooth initial data 
\[
\phi(0,x)=\phi_0(x),\qquad \frac{\partial\phi}{\partial t}(0,x)=\phi_1(x)
\]
has a smooth unique solution.
\end{Satz}

To prove Theorem \ref{satz-cauchy-magnetic} we adapt the proof
obtained for (standard) wave maps from two-dimensional Minkowski space \cite{MR596432}.

We will divide the proof of Theorem \ref{satz-cauchy-magnetic} into several steps.

\begin{Lem}
The Cauchy problem \eqref{cauchy-magnetic} is equivalent to the system
\begin{align}
\label{system-uv}
\frac{\partial u^i}{\partial\eta}+(\Gamma^i_{jk}+\frac{1}{2}Z^i(\partial_{y^k}\wedge\partial_{y^j}))u^jv^k&=0,\qquad \frac{\partial y^i}{\partial\xi}=u^i, \qquad i=1,\ldots n,\\
\nonumber\frac{\partial v^i}{\partial\xi}+(\Gamma^i_{jk}+\frac{1}{2}Z^i(\partial_{y^k}\wedge\partial_{y^j}))v^ju^k&=0,\qquad \frac{\partial z^i}{\partial\eta}=v^i,\qquad i=1,\ldots n
\end{align}
for functions \(u,v,y,z\colon\R^{1,1}\to\R^n\) with the initial conditions
\begin{align}
\label{lightcone-b-initial}
y^i(0,x)=z^i(0,x)=\phi^i_0(x),\qquad
u^i(0,x)=\frac{\partial\phi_0^i}{\partial x}(x)+\phi_1^i(x),\qquad  v^i(0,x)=-\frac{\partial\phi_0^i}{\partial x}(x)+\phi_1^i(x).
\end{align}
\end{Lem}
\begin{proof}
This follows by a direct calculation using lightcone coordinates.
\end{proof} 

In the following let \(U\) be a coordinate system in the target \(N\) with coordinates \(y^i,i=1\ldots n\), and \(\tilde U\)
a subset of \(U\) fixed by the requirement
\[
\sum_{i=1}^n(y^i)^2\leq 1.
\]
Let \(a\) be a positive constant such that
\[
g(y)(\lambda,\lambda)\geq a^2|\lambda|^2
\]
holds for all \(y\in\tilde U\) and all \(\lambda\in\R^n\).

We set 
\begin{align*}
M:=\sup_{|x|\leq L}\big\{\sqrt{g(\phi_0)(dy(\partial_\xi),dy(\partial_\xi))},\sqrt{g(\phi_0)(dy(\partial_\eta),dy(\partial_\eta))}\big\}
\end{align*}
for a large number \(L\).

As a next step we prove the existence of a solution to the system \eqref{system-uv} on the characteristic triangle \(\Lambda_k\),
which we define by
\[
-\frac{1}{2}k\leq-\eta\leq\xi\leq \frac{1}{2}k
\]
with 
\[
k:=\min\{L,\frac{a}{2M\sqrt{n}}\}.
\]

\begin{Prop}
There exists a smooth solution to the Cauchy problem \eqref{cauchy-magnetic} on the characteristic triangle \(\Lambda_k\).
\end{Prop}
\begin{proof}
To establish the claim we will make use of the system \eqref{system-uv}.
Let \(u_0^i,v_0^i\) be smooth functions satisfying the initial condition \eqref{lightcone-b-initial} 
\[
\frac{\partial y^i_0}{\partial\xi}=u_0^i,\qquad \frac{\partial z_0^i}{\partial\eta}=v_0^i,\qquad y_0^i(0,x)=z_0^i(0,x)=\phi^i_0(x)
\]
and the restrictions
\begin{align}
\label{constraint}
\sum_{i=1}^n(y_0^i)^2\leq 1, \qquad \sum_{i=1}^n(z_0^i)^2\leq 1.
\end{align}
Suppose we have constructed \(y^i_{m-1},z^i_{m-1},u^i_{m-1},v^i_{m-1}\) that satisfy the initial condition \eqref{lightcone-b-initial}.
Moreover, suppose that \(y^i_{m-1}\) and \(z^i_{m-1}\) satisfy the constraint \eqref{constraint}.
We define \(y^i_{m},z^i_{m},u^i_{m},v^i_{m}\) by the equations
\begin{align}
\label{lightcone-b-iteration}
\frac{\partial u_m^i}{\partial\eta}+(\Gamma^i_{jk}(z_{m-1})+\frac{1}{2}Z^i(z_{m-1})(\partial_{y^k}\wedge\partial_{y^j}))u_m^jv_{m-1}^k=0,
\qquad\frac{\partial y_m^i}{\partial\xi}=u_m^i, \\
\nonumber\frac{\partial v_m^i}{\partial\xi}+(\Gamma^i_{jk}(y_{m-1})
+\frac{1}{2}Z^i(y_{m-1})(\partial_{y^k}\wedge\partial_{y^j}))v_m^ju_{m-1}^k=0,\qquad\frac{\partial z_m^i}{\partial\eta}=v_m^i.
\end{align}
Note that the equations \eqref{lightcone-b-iteration} are linear. Hence, their solutions exist in the whole characteristic
triangle \(\Lambda_k\).

In addition, the equations \eqref{lightcone-b-iteration} also have a nice geometric interpretation.
We can think of these equations describing two parallel vector fields \(u_m^i,v_m^i\)
along the curves \(z^i=z^i_{m-1}(\xi_0,\eta)\) and \(y^i=y^i_{m-1}(\xi,\eta_0)\) for a
metric connection with totally antisymmetric torsion.

Since parallel translation is an isometry (this is not affected by the presence of torsion),
we have
\begin{align*}
g(z_{m-1})(u_m,u_m)=g(z_{m-1})(u_m,u_m)\big|_{t=0},\qquad g(z_{m-1})(v_m,v_m)=g(z_{m-1})(v_m,v_m)\big|_{t=0}.
\end{align*}
Thus, we obtain the estimates
\begin{align*}
|u^i_m|\leq\frac{M}{a},\qquad |v^i_m|\leq\frac{M}{a},\qquad i=1,\ldots,n.
\end{align*}
Consequently, we also find
\begin{align*}
\sum_{i=1}^n(y_m^i)^2\leq 4n\frac{M^2}{a^2}k^2\leq 1,\qquad \sum_{i=1}^n(z_m^i)^2\leq 1.
\end{align*}
We can deduce that \(y^i_m,z^i_m\) remain in the set \(\tilde{U}\) and that \(u_m,v_m\) are bounded.
Consequently, the sequences \(\{u_m\},\{v_m\},\{z_m\},\{y_m\}\) and the sequences of their partial derivatives
converge uniformly on the characteristic triangle \(\Lambda_k\). Hence, the limit 
\[
y(t,x)=\lim_{m\to\infty}y_m(t,x)
\]
is a solution to the Cauchy problem \eqref{cauchy-magnetic} on \(\Lambda_k\).
\end{proof}

\begin{Lem}[Global Existence]
The solution constructed on the characteristic triangle \(\Lambda_k\) can be extended to
the whole Minkowski space.
\end{Lem}

\begin{proof}
To prove that there exists a global smooth solution of \eqref{cauchy-magnetic} it suffices to know that
there exists a smooth solution on \(\Lambda_k\). This can be seen as follows:
We know that the Cauchy problem \eqref{cauchy-magnetic} is solvable up to \(t=k\). We perform the following rescaling of the initial data
\[
\phi_0(k,x)=\phi(\sigma k,x),\qquad \phi_1(k,x)=\frac{\partial\phi}{\partial t}(\sigma k,x),\qquad -L+\sigma k\leq x\leq L-\sigma k,
\]
where \(\sigma\geq 1\) is a fixed positive integer, defined by \(t=\sigma k\).
We infer that the Cauchy problem \eqref{cauchy-magnetic} is solvable up to \(t=(\sigma+1)k\).
By iteration of this procedure we conclude that there exists a global smooth solution
of \eqref{cauchy-magnetic}. 
\end{proof}

In the following we will need the extrinsic version of \eqref{euler-lagrange-magnetic}, which is given by
\begin{equation}
\label{extrinsic-magnetic}
\frac{\partial^2u}{\partial t^2}-\frac{\partial^2u}{\partial x^2}=\sff(du,du)+Z(du(\partial_t)\wedge du(\partial_x)),
\end{equation}
where \(u\colon\R^{1,1}\to\R^q\).

\begin{Lem}[Uniqueness]
The global solution to the Cauchy problem \eqref{cauchy-magnetic} constructed in the last Lemma is unique.
\end{Lem}
\begin{proof}
Suppose we have two solutions \(u,v\) of the extrinsic equation for wave maps with two-form potential \eqref{extrinsic-magnetic}
that have the same initial data, that is \(u=v\) and \(du=dv\) at \(t=0\).
We set \(w:=u-v\). Making use of the extrinsic version for wave maps with two-form potential we
find by a direct calculation
\begin{align*}
\frac{d}{dt}\int_\R(\big|\frac{\partial w}{\partial t}\big|^2+\big|\frac{\partial w}{\partial x}\big|^2)dx
=&\int_\R\langle\frac{\partial w}{\partial t},\sff(u)(du,du)-\sff(v)(dv,dv)\rangle dx\\
&+\int_\R\langle\frac{\partial w}{\partial t},Z(u)(du(\partial_t)\wedge du(\partial_x))-Z(v)(dv(\partial_t)\wedge dv(\partial_x))\rangle dx.
\end{align*}
We have a pointwise bound on \(du\) and \(dv\) from \eqref{lightcone-b-iteration}, which allows us to estimate
\begin{align*}
|\langle \sff(u)(du,du)-\sff(v)(dv,dv),\frac{\partial w}{\partial t}\rangle|&\leq C(|dw|^2+|dw||w|), \\
|\langle Z(u)(du(\partial_t)\wedge du(\partial_x))-Z(v)(dv(\partial_t)\wedge dv(\partial_x)),\frac{\partial w}{\partial t}\rangle|
&\leq C(|dw|^2+|dw||w|).
\end{align*}
Thus, we find
\begin{align*}
\frac{d}{dt}\int_\R(w^2+|dw|^2)dx\leq C \int_\R(w^2+|dw|^2)dx.
\end{align*}
By integrating the differential equation and using that \(w=dw=0\) at \(t=0\) we obtain the claim.
\end{proof}

This finishes the proof of Theorem \ref{satz-cauchy-magnetic}.

\begin{Bem}
The action functional for the full bosonic string  presented in the introduction
contains an additional scalar potential. If we try to extend the strategy of proof
used in this section to the latter it turns out that we do not 
succeed due to the presence of the scalar potential.

On the one hand we lose the scale invariance of the Euler-Lagrange
equations and secondly, we cannot rewrite the Euler-Lagrange equation into a system which
has a geometric interpretation as parallel transport. 
\end{Bem}
\par\medskip
The proof of Theorem \ref{satz-cauchy-magnetic} presented above is of a geometric nature.
In the following we want to sketch briefly how Theorem \ref{satz-cauchy-magnetic} can also be proven making use of analytical methods, which
is inspired from \cite{MR933231}. These methods have the advantage that they can also be applied for non-smooth initial data.

For a solution \(u\) of \eqref{extrinsic-magnetic} we define the energy
\[
e(t,x):=\frac{1}{2}\big|\frac{\partial u}{\partial t}\big|^2+\frac{1}{2}\big|\frac{\partial u}{\partial x}\big|^2.
\]
If we test \eqref{extrinsic-magnetic} with \(\frac{\partial u}{\partial t}\) we find
\begin{align*}
0=&\langle\sff(du,du),\frac{\partial u}{\partial t}\rangle+\langle Z(\frac{\partial u}{\partial t}\wedge\frac{\partial u}{\partial x}),\frac{\partial u}{\partial t}\rangle \\
\nonumber=&\langle\frac{\partial^2u}{\partial t^2},\frac{\partial u}{\partial t}\rangle-\langle\frac{\partial^2u}{\partial x^2},\frac{\partial u}{\partial t}\rangle \\
\nonumber=&\frac{1}{2}\frac{\partial}{\partial t}\big|\frac{\partial u}{\partial t}\big|^2+\frac{1}{2}\frac{\partial}{\partial t}\big|\frac{\partial u}{\partial x}\big|^2
-\frac{\partial}{\partial x}\langle\frac{\partial u}{\partial x},\frac{\partial u}{\partial t}\rangle,
\end{align*}
which yields the identity
\begin{align}
\label{conservation-a}
\frac{\partial}{\partial t}e(t,x)=\frac{\partial}{\partial x}\langle\frac{\partial u}{\partial x},\frac{\partial u}{\partial t}\rangle.
\end{align}
Interchanging the roles of \(t\) and \(x\) a similar calculation leads to
\begin{align}
\label{conservation-b}
\frac{\partial}{\partial x}e(t,x)=\frac{\partial}{\partial t}\langle\frac{\partial u}{\partial x},\frac{\partial u}{\partial t}\rangle.
\end{align}
Differentiating \eqref{conservation-a} with respect to \(t\) and \eqref{conservation-b} with respect to \(x\) yields the identity
\begin{align}
\label{box-energy-magnetic}
\big(\frac{\partial^2}{\partial t^2}-\frac{\partial^2}{\partial x^2}\big) e(t,x)=0,
\end{align}
which means that the energy \(e(t,x)\) satisfies the linear wave equation in two-dimensional Minkowski space.
Consequently, we get the following energy equality
\begin{align*}
\tilde{E}(t):=\frac{1}{2}\int_\R\big(\big|\frac{\partial}{\partial t} e(t,x)\big|^2+|\frac{\partial}{\partial x} e(t,x)\big|^2\big)dx=\tilde{E}(0).
\end{align*}
Hence, by the Sobolev embedding theorem, the energy \(e(t,x)\) is uniformly bounded in terms of initial data. Depending on the regularity
of the initial data we get a global solution of the corresponding regularity of the Cauchy problem \eqref{cauchy-magnetic}.

\section{A global weak solution for wave maps with scalar and two-form potential}
In this section we derive an existence result for wave maps with scalar and two-form potential. 
In order to take into account the additional scalar potential we have to make use of a more general approach than in the last section. 

We use the notation \(\partial_t\) and \(\partial_x\) to represent partial derivatives with respect to \(t\) and \(x\).
We set \(\Box:=\partial_t^2-\partial_x^2\).
Recall the extrinsic version of the equation for wave maps with scalar and two-form potential
\begin{align}
\label{boxu}
\Box u=\sff(du,du)+Z(\partial_t u\wedge\partial_x u)+\nabla V(u)
\end{align}
for a function \(u\colon\R^{1,1}\to N\subset \R^q\). In this section we study the well-posedness of the Cauchy problem
with initial data 
\begin{align*}
(u,\partial_tu)|_{t=0}=(u_0,u_1),
\end{align*}
where \(u_0(x)\in N\subset\R^q\) and \(u_1(x)\in T_{u_0(x)N}\subset\R^q\).
It turns out that we have to consider a bigger class of initial data compared to the last section.

\begin{Bem}
For an arbitrary scalar potential \(V(u)\) solutions of \eqref{boxu} will not be invariant under scaling.
However, in the special case that we perform the rescaling
\[
u(t,x)\to u_\lambda(t,x):=\lambda^\beta u(\lambda t,\lambda t),\qquad \lambda>0, \qquad \beta\in\R
\]
and the potential \(V(u)\) is such that it satisfies
\[
\nabla V(u)\to \lambda^{2+\beta}\nabla V(u_\lambda),
\]
then solutions of \eqref{boxu} are scale invariant.
\end{Bem}

\subsection{Existence of a global weak solution}
By controlling the \(H^2\)-norm of \(u\) we will derive the following existence result:

\begin{Satz}
\label{theorem-full}
Let \((N,g)\) be a compact Riemannian manifold.
Suppose that the scalar potential \(V\) is of class \(C^2(N,\R)\).
For any initial data of the regularity
\begin{align*}
(u_0,u_1)\in H^2\times H^1(\R,TN)
\end{align*}
the Cauchy problem \eqref{boxu} admits a unique solution of class \(H^2\) for all times \(t\).
\end{Satz}

First of all, we note that for a continuous potential \(V(u)\) and due to the assumption that \(N\) is compact,
we have the following estimate
\begin{align*}
\int_\R(\big|\frac{\partial u}{\partial t}\big|^2+\big|\frac{\partial u}{\partial x}\big|^2)dx=\int_\R|du|^2dx=E(0)+2\int_\R V(u)dx\leq C.
\end{align*}

For a general treatment of semilinear wave equations of the form \eqref{boxu} we refer to \cite[Chapter 5]{MR1674843}.
Moreover, by \cite[Remark 5.3]{MR1674843} there exists a unique local (in time) solution to the Cauchy problem \eqref{boxu}.

We set
\[
E_2(t):=\frac{1}{2}\int_\R(|\partial^2_xu|^2+2|\partial_x\partial_tu|^2+|\partial^2_tu|^2)dx.
\]

To extend the local solution to a global one we establish the following
\begin{Prop}
Let \(u\colon[0,T)\times\R\to N\subset \R^q\) be a weak solution of \eqref{boxu}.
Then the following inequality holds
\begin{align}
\label{evolution-boxu-h2}
\frac{dE_2(t)}{dt}\leq C\big(E_2(t)\int_\R|du|^2dx+E_2(t)^\frac{1}{2}\big(\int_\R|du|^2dx\big)^\frac{1}{2}\big),
\end{align}
where the constant \(C\) depends on \(Z,\nabla Z,\nabla V\) and \(\hess V\).
\end{Prop}
\begin{proof}
We calculate 
\begin{align*}
\frac{dE_2(t)}{dt}=\int_\R(\langle\partial_x(\Box u),\partial_x\partial_tu\rangle+\langle\partial_t(\Box u),\partial^2_tu\rangle)dx.
\end{align*}
Using \eqref{boxu} we get
\begin{align*}
\langle\partial_x(\Box u),\partial_x\partial_tu\rangle=\langle\partial_x(\sff(du,du)),\partial_x\partial_tu\rangle
+\langle\partial_x(Z(\partial_tu\wedge\partial_xu)),\partial_x\partial_tu\rangle
+\langle\partial_x(\nabla V(u)),\partial_x\partial_tu\rangle
\end{align*}
and similarly for the derivatives with respect to \(t\).
First of all, we note that
\begin{align*}
\langle\partial_x(\sff(du,du)),\partial_x\partial_tu\rangle=&\langle(\partial_{du(\partial_x)}\sff)(du,du),\partial_x\partial_tu\rangle
+2\langle\sff(d\partial_xu,du),\partial_x\partial_tu\rangle \\
=&\langle(\partial_{du(\partial_x)}\sff)(du,du),\partial_x\partial_tu\rangle-2\langle(\partial_{du(\partial_x)}\sff)(d\partial_xu,du),\partial_tu\rangle,
\end{align*}
where we used that \(\partial_t\perp\sff\). Hence, we get the estimate
\begin{align*}
|\langle\partial_x(\sff(du,du)),\partial_x\partial_tu\rangle+\langle\partial_t(\sff(du,du)),\partial^2_tu\rangle|\leq C|du|^3|d^2u|.
\end{align*}
Moreover, we have
\begin{align*}
\langle\partial_x(Z(\partial_tu\wedge\partial_xu))&,\partial_x\partial_tu\rangle+\langle\partial_t(Z(\partial_tu\wedge\partial_xu)),\partial^2_tu\rangle \\
=&\langle(\partial_{du(\partial_x)}Z)(\partial_tu\wedge\partial_xu),\partial_x\partial_tu\rangle+\Omega(\partial_tu,\partial^2_xu,\partial_x\partial_tu) \\
&+\langle(\partial_{du(\partial_t)}Z)(\partial_tu\wedge\partial_xu),\partial^2_tu\rangle+\Omega(\partial_tu,\partial_x\partial_tu,\partial^2_tu) \\
=&\langle(\partial_{du(\partial_x)}Z)(\partial_tu\wedge\partial_xu),\partial_x\partial_tu\rangle+\langle(\partial_{du(\partial_t)}Z)(\partial_tu\wedge\partial_xu),\partial^2_tu\rangle \\
&+\Omega(\partial_tu,\partial_x\partial_tu,\Box u).
\end{align*}
Hence, we may conclude
\begin{align*}
|\langle\partial_x(Z(\partial_tu\wedge\partial_xu)),\partial_x\partial_tu\rangle+\langle\partial_x(Z(\partial_tu\wedge\partial_xu)),\partial^2_tu\rangle|\leq C(|du|^3|d^2u|+|du||d^2u|),
\end{align*}
where we used \eqref{boxu} again. Finally, we note that
\begin{align*}
\langle\partial_x(\nabla V(u)),\partial_x\partial_tu\rangle=\hess V(u)(\partial_xu,\partial_x\partial_tu), \qquad
\langle\partial_t(\nabla V(u)),\partial^2_tu\rangle=\hess V(u)(\partial_tu,\partial^2_tu).
\end{align*}
We estimate
\begin{align*}
\big\||du|^3|d^2u|\big\|_{L^1(\R)}\leq & \big\||du|^3\big\|_{L^2(\R)}\big\||d^2u|\big\|_{L^2(\R)} \\
\leq &C\big\|d|du|^3\big\|_{L^\frac{2}{3}(\R)}\big\||d^2u|\big\|_{L^2(\R)} \\
\leq &C\big\||du|^2|d^2u|\big\|_{L^\frac{2}{3}(\R)}\big\||d^2u|\big\|_{L^2(\R)} \\
\leq &C\big\||du|\big\|^2_{L^2(\R)}\big\||d^2u|\big\|^2_{L^2(\R)},
\end{align*}
where we used the Sobolev embedding theorem in one dimension in the second step.
Thus, we obtain
\begin{align*}
\frac{dE_2(t)}{dt}\leq C\big(\int_\R|du|^2dx\int_\R|d^2u|^2dx+\int_\R|du||d^2u|dx\big),
\end{align*}
which completes the proof.
\end{proof}

\begin{Cor}
Let \(u\colon[0,T)\times\R\to N\subset \R^q\) be a weak solution of \eqref{boxu}.
Then the following inequality holds
\begin{align*}
E_2(t)\leq C_1e^{C_2T},
\end{align*}
where the positive constants \(C_1\) and \(C_2\) depend on \(Z,\nabla Z,V,\nabla V,\hess V\) and the initial data.
\end{Cor}
\begin{proof}
Using the conserved energy \(\int_\R|du|^2dx\) and dividing by \(E_2(t)^\frac{1}{2}\) equation \eqref{evolution-boxu-h2} becomes an inequality of the form
\[
\frac{d}{dt}f(t)\leq c_1f(t)+c_2,
\]
where \(c_1\) and \(c_2\) are positive constants.
The claim follows by application of the Gronwall Lemma.
\end{proof}
Having gained control over the spatial \(H^2\)-norm of \(u\) we are able to show
the global existence of a solution to \eqref{boxu}.

\begin{Cor}
Let \(u\colon[0,T)\times\R\to N\subset \R^q\) be a weak solution of \eqref{boxu}.
The solution exists for all \(t\in[0,\infty)\).
\end{Cor}
\begin{proof}
Since we have a bound on the spatial \(H^2\)-norm of \(u\) we can extend the solution 
for all \(t\in[0,\infty)\).
\end{proof}

Finally, we address the issue of uniqueness.

\begin{Prop}[Uniqueness]
Suppose \(u,v\colon\R^{1,1}\to N\subset \R^q\) are two \(H^2\) solutions of \eqref{boxu}. If they coincide at \(t=0\), then they coincide for all times.
\end{Prop}
\begin{proof}
We set \(w:=u-v\). By assumption both \(u\) and \(v\) are solutions of \eqref{boxu}. By the last Proposition
we know that the \(H^2\)-norm of both \(u\) and \(v\) is finite. Making use of the Sobolev embedding \(W^{2,2}(\R)\subset C^{1,\frac{1}{2}}(\R)\)
we obtain a pointwise bound on \(|du|\) and \(|dv|\).
We will use the notation \(u_x\) to denote the partial derivative of \(u\) with respect to \(x\) and similarly for 
the derivative with respect to \(t\).
Performing a direct calculation we obtain
\begin{align*}
\frac{d}{dt}\int_\R|dw|^2dx=\int_\R&\big(\langle w_t,\sff(u)(du,du)-\sff(v)(dv,dv)\rangle\\
&+\langle w_t,Z(u)(u_t\wedge u_x)-Z(v)(v_t\wedge v_x)\rangle 
+\langle w_t,\nabla V(u)-\nabla V(v)\rangle\big)dx\\
&:=A_1+A_2+A_3.
\end{align*}
We rewrite the \(A_1\) term as follows
\[
\sff(u)(du,du)-\sff(v)(dv,dv)=\sff(u)(du,du)-\sff(v)(du,du)+\sff(v)(du+dv,dw).
\]
This yields
\[
\langle w_t,(\sff(u)-\sff(v))(du,du)\rangle\leq C|du|^2|w||dw|.
\]
Using the orthogonality \(\partial_t\perp\sff\) we find
\[
\langle w_t,\sff(v)(du,dw)\rangle=\langle u_t,\sff(v)(du,dw)\rangle=\langle u_t,(\sff(v)-\sff(u))(du,dw)\rangle\leq C|du|^2|w||dw|.
\]
The same argument also holds for the term involving \(dv\).
Consequently, the \(A_1\)-term can be bounded as 
\[
A_1\leq C\int_\R|du|^2|w||dw|dx\leq C\int_\R|w||dw|dx,
\]
where we used the pointwise bound on \(|du|\).
To control the \(A_2\)-term we rewrite
\begin{align*}
\langle w_t,Z(u)(u_t\wedge u_x)-Z(v)(v_t\wedge v_x)\rangle=
&\langle w_t,(Z(u)-Z(v))(u_t\wedge u_x)\rangle 
+\langle w_t,Z(v)(w_t\wedge u_x)\rangle \\
&+\langle w_t,Z(v)(v_t\wedge w_x)\rangle
\end{align*}
such that we can estimate
\[
B\leq C\int_\R(|w||dw||du|^2+|dw|^2(|du|+|dv|))dx.
\]
By assumption \(\nabla V\) is continuous, hence we may estimate
\begin{align*}
|\langle w_t,\nabla V(u)-\nabla V(v)\rangle|\leq C|dw||w|.
\end{align*}
Applying Young's inequality several times we are led to
\begin{align*}
\frac{d}{dt}\int_\R(|w|^2+|dw|^2)dx\leq C\int_\R(|w|^2+|dw|^2)dx,
\end{align*}
which can be integrated as
\begin{align*}
\int_\R(|u-v|^2+|du-dv|^2)dx\big|_{t=T}\leq e^{CT}\int_\R(|u(0)-v(0)|^2+|du(0)-dv(0)|^2)dx.
\end{align*}
By assumption we have \(u(0)=v(0)\) and \(du(0)=dv(0)\). Consequently,
we have \(u(t)=v(t)\) for all \(t\in [0,\infty)\).
\end{proof}

This finishes the proof of Theorem \ref{theorem-full}.

\begin{Bem}
The wave map equation has been studied for less regular initial data,
see \cite{MR1759884} and \cite{MR1663216}. However, a sophisticated analysis performed in 
these references shows that the Cauchy problem is ill-posed then.
\end{Bem}

\subsection{A priori estimates}
In this subsection we derive an a priori estimate for weak solutions of \eqref{boxu},
which heavily makes use of the special structure of the right hand side of \eqref{boxu}.
Our result is inspired by a similar estimate for wave maps from \cite{MR1381973}.

We define the following quantity
\begin{align}
\label{definition-F}
F_\pm(u):=(\partial_tu\pm\partial_xu),
\end{align}
which leads to the following
\begin{Prop}
Let \(u\colon [0,T)\times\R\to N\) be a weak solution of \eqref{boxu}.
Then the following energy inequality holds
\begin{align}
\int_0^T\int_\R (F^2_-(u)F^2_+(u)-4V^2(u))dxdt 
\leq 2 E(0)^2+2\int_{x\geq y}\int (V(u)(|\partial_xu|^2+|\partial_tu|^2))dxdy\big|_{t=0}.
\end{align}
\end{Prop}

\begin{proof}
We rewrite the extrinsic version of the Euler-Lagrange equation \eqref{boxu} as
\begin{align*}
(\partial_t+\partial_x)(\partial_t-\partial_x)u=\sff(du,du)+Z(\partial_tu\wedge\partial_xu)+\nabla V(u).
\end{align*}
Taking the scalar product with \((\partial_t\pm\partial_x) u\) we find
\begin{align*}
(\partial_t\pm\partial_x)\big|(\partial_t\mp\partial_x)u\big|^2=2(\partial_t\mp\partial_x)V(u).
\end{align*}
Making use of \eqref{definition-F} we can rewrite this as a conservation law in the form
\[
\partial_t(F^2_\pm(u)-2V(u))=\pm\partial_x(F^2_\pm(u)+2V(u)).
\]
Note that
\[
\int_\R F^2_\pm(u)dx=\int_\R(|\partial_xu|^2+|\partial_tu|^2\pm 2\langle \partial_xu,\partial_tu\rangle)dx\leq 2 E(0)+2\int_\R V(u)dx.
\]

We define the quantity
\[
Z(t):=\int_{x\geq y}\int (F^2_-(u)-2V(u))(F^2_+(u)-2V(u))dxdy,
\]
which satisfies the bound 
\[
Z(t)\leq\int_\R (F^2_-(u)-2V(u))dx\int_\R (F^2_+(u)-2V(u))dx\leq 4E(0)^2.
\]
We compute
\begin{align*}
\frac{d}{dt}Z(t)=&\int_{x\geq y}\int\big(\partial_t(F^2_-(u)-2V(u))(F^2_+(u)-2V(u)) \\
&+(F^2_-(u)-2V(u))\partial_t (F^2_+(u)-2V(u))\big)dxdy \\
=&-\int_{-\infty}^\infty\big(\int_y^\infty\partial_x(F^2_-(u)+2V)dx\big)(F^2_+(u)-2V(u))dy \\
&+\int_{-\infty}^\infty\big(\int^x_\infty\partial_y(F^2_+(u)+2V(u))dy\big)(F^2_-(u)-2V(u))dx\\
=&2\int_\R (F^2_-(u)F^2_+(u)-4V^2(u))dx.
\end{align*}
The result follows by integration with respect to \(t\). 
\end{proof}

\begin{Bem}
The last Proposition gives a \(H^2\) bound of \(u\) in spacetime in terms of the potential \(V(u)\) and the initial data.
\end{Bem}

\begin{Bem}
If we perform the same calculation that led to \eqref{box-energy-magnetic} taking into account the
scalar potential, we find
\begin{align}
\Box e(t,x)=(\frac{\partial^2}{\partial x^2}+\frac{\partial^2}{\partial t^2})V(u), 
\end{align}
where
\[
e(t,x):=\frac{1}{2}\big|\frac{\partial u}{\partial t}\big|^2+\frac{1}{2}\big|\frac{\partial u}{\partial x}\big|^2.
\]
Making use of Duhamel's Principle we can write down a formal solution to this equation, which would
give a pointwise bound on \(e(t,x)\). However, it seems that these bounds cannot be used to improve
Theorem \ref{theorem-full}.
\end{Bem}

\par\medskip
\textbf{Acknowledgements:}
The author gratefully acknowledges the support of the Austrian Science Fund (FWF) 
through the START-Project Y963-N35 of Michael Eichmair.

\bibliographystyle{plain}
\bibliography{mybib}
\end{document}